\input amstex
\input amsppt.sty   
\input epsf.tex
\magnification=\magstep1
\hsize 30pc
\vsize 47pc
\def\nmb#1#2{\ifx#1!{\rm (#2)}\else#2\fi} 
\def\cit#1#2{\ifx#1!\cite{#2}\else#2\fi} 
\def\ign#1{}             
\redefine\o{\circ}
\define\X{\frak X}

\define\ep{\varepsilon}

\redefine\i{^{-1}}
\define\x{\times}
\define\Fl{\operatorname{Fl}}
\let\on\operatorname
\def\today{\ifcase\month\or
 January\or February\or March\or April\or May\or June\or
 July\or August\or September\or October\or November\or December\fi
 \space\number\day, \number\year}
\topmatter
\title  The flow completion of a manifold with vector field
\endtitle
\author  Franz W. Kamber, Peter W. Michor  \endauthor
\affil
Erwin Schr\"odinger International Institute of Mathematical Physics, 
Wien, Austria
\endaffil
\address
Franz W. Kamber: 
Department of Mathematics, 
University of Illinois, 
1409 West Green Street, 
Urbana, IL 61801, USA
\endaddress
\email kamber\@math.uiuc.edu \endemail
\address
P\. W\. Michor: Institut f\"ur Mathematik, Universit\"at Wien,
Strudlhofgasse 4, A-1090 Wien, Austria; {\it and:} 
Erwin Schr\"odinger Institut f\"ur Mathematische Physik,
Boltzmanngasse 9, A-1090 Wien, Austria
\endaddress
\email michor\@pap.univie.ac.at \endemail
\date {\today} \enddate
\thanks Supported by 
Erwin Schr\"odinger International Institute of Mathematical Physics, 
Wien, Austria.
FWK was supported in part by The National Science Foundation under 
Grant No. DMS-9504084.
PWM was supported  
by `Fonds zur F\"orderung der wissenschaftlichen  
Forschung, Projekt P~14195~MAT'.
\endthanks
\abstract For a vector field $X$ on a smooth manifold $M$ there 
exists a smooth but not necessarily Hausdorff manifold $M_\Bbb R$ and 
a complete vector field $X_\Bbb R$ on it which is the universal 
completion of $(M,X)$. 
\endabstract
\endtopmatter

\document

\proclaim{\nmb.{1}. Theorem} 
Let $X\in\X(M)$ be a smooth vector field on a 
(connected)
smooth manifold $M$. 

Then there exists 
a universal flow completion $j:(M,X)\to (M_\Bbb R,X_\Bbb R)$ of $(M,X)$. 
Namely, there exists
a (connected) smooth not necessarily Hausdorff manifold $M_\Bbb R$, 
a complete vector field $X_\Bbb R\in \X(M_\Bbb R)$, and an embedding 
$j:M\to M_\Bbb R$ onto an open submanifold such that $X$ and 
$X_\Bbb R$ are $j$-related: $Tj\o X=X_\Bbb R\o j$. Moreover, for any 
other equivariant morphism $f:(M,X)\to (N,Y)$ for a manifold $N$ and 
a complete vector field $Y\in X(N)$ there exists a unique equivariant 
morphism $f_\Bbb R:(M_\Bbb R,x_\Bbb R)\to (N,Y)$ with 
$f_\Bbb R\o j=f$. The leave spaces $M/X$ and $M_\Bbb R/X_\Bbb R$ are 
homeomorphic. 
\endproclaim

\demo{Proof}
Consider the manifold $\Bbb R\x M$ with coordinate function $s$ on 
$\Bbb R$, the vector field $\bar X:=\partial_s\x X\in\X(\Bbb R\x M)$, 
and let $M_\Bbb R:= \Bbb R\x_{\bar X}M$ be the orbit space (or leaf 
space) of the vector field $\bar X$. 

Consider the flow mapping 
$\Fl^{\bar X}:\Cal D(\bar X)\to \Bbb R\x M$, given by 
$\Fl^{\bar X}_t(s,x)=(s+t,\Fl^X_t(x))$, where the domain of definition 
$\Cal D(\bar X)\subset \Bbb R\x (\Bbb R\x M)$ is an open 
neighbourhood of $\{0\}\x(\Bbb R\x M)$ with the property that 
$\Bbb R\x \{x\}\cap \Cal D(\bar X)$ is an open interval times 
$\{x\}$. 

For each $s\in \Bbb R$ we consider the mapping 
$$
j_s:M @>{\on{ins_t}}>> \{s\}\x M\subset \Bbb R\x M @>{\pi}>> 
\Bbb R\x _{\bar X}M= M_\Bbb R.
$$
Each mapping $j_s$ is injective: A trajectory of $\bar X$ can 
meet $\{s\}\x M$ at most once since it projects onto the unit speed 
flow on $\Bbb R$. 

Obviously, the image $j_s(M)$ is open in $M_\Bbb R$ in the 
quotient topology:
If a trajectory hits $\{s\}\x M$ in a point $(s,x)$, let $U$ be an open 
neighborhood of $x$ in $M$ such that 
$(-\ep,\ep)\x(s-\ep,s+\ep)\x U\subset \Cal D(\bar X)$. Then the 
trajectories hitting $(s-\ep,s+\ep)\x U$ fill a flow invariant open 
neighborhood which projects on an open neighborhood of $j_s(x)$ in 
$M_\Bbb R$ which lies in $j_s(M)$. This argument also shows that 
$j_s$ is a homeomorphism onto its image in $M_\Bbb R$.

Let us use the mappings $j_s:M\to M_\Bbb R$ as charts. The chart 
change then looks as follows: For $r<s$ the set 
$(j_s)\i(j_r(M))\subset M$ is just the open subset of all
$x\in M$ such that $[0,s-r]\x\{(s,x)\}\subset \Cal D(\bar X)$, and 
$(j_s)\i\o j_r$ is given by $\Fl^X_{s-r}$ on this set. Thus the chart 
changes are smooth. 

Consider the flow $(t,(s,x))\mapsto (s+t,x)$ on 
$\Bbb R\x M$ which commutes with the flow of $\bar X$ and thus 
induces a flow on the leave space $M_\Bbb R=\Bbb R\x_{\bar X}M$. 
Differentiating this flow we get a vector field $X_\Bbb R$ on 
$M_\Bbb R$. 

The construction $(M,X)\mapsto (M_\Bbb R,X_\Bbb R)$ is a functor from 
the category of smooth Hausdorff manifolds with vector-fields and 
smooth mappings intertwining the vector fields into the category 
of possibly non-Hausdorff manifolds with complete smooth vector 
fields and smooth mappings intertwining these fields. 
For a pair $(M,X)$ with $X$ a complete vector field the flow 
completion $(M_\Bbb R,X_\Bbb R)$ is equivariantly diffeomorphic to 
$(M,X)$ since then any of the charts $j_s:M\to M_\Bbb R$ is also 
surjective. From this the universal property follows. 
\qed\enddemo

\subhead\nmb.{2}. Example \endsubhead
Let $(M,X)=(\Bbb R^2\setminus\{0\},\partial_x)$.
The trajectories of $X$ on $M$ and of $\bar X$ on $\Bbb R\x M$ in the 
slices $y=\text{ constant}$ for $y=0$ and $y\ne0$ then look as follows:
$$\alignat3
&M &\qquad &\Bbb R\x M,\quad y=0  &\qquad &\Bbb R\x M,\quad y\ne 0 \\
&{\epsfxsize=2cm\epsfbox{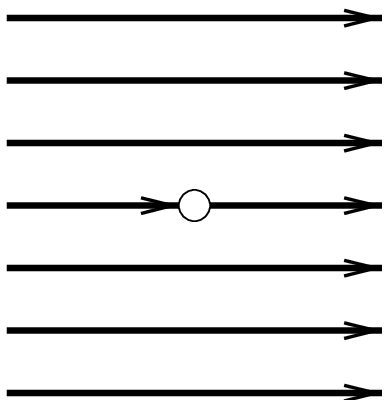}} &\qquad
&{\epsfxsize=2cm\epsfbox{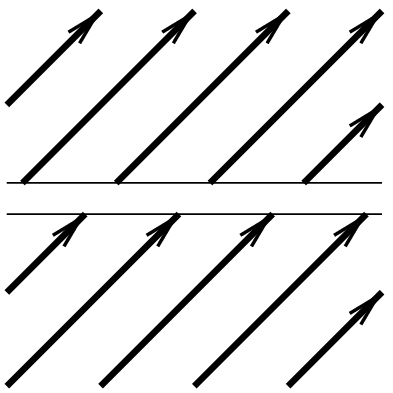}} &\qquad
&{\epsfxsize=2cm\epsfbox{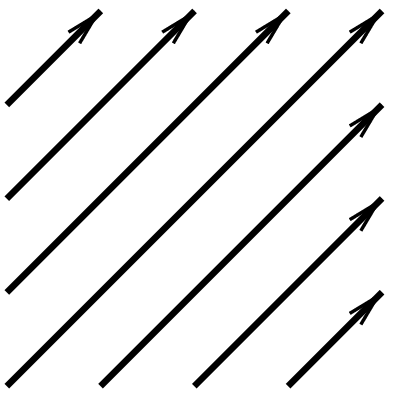}} 
\endalignat$$ 
The smooth manifold $M_\Bbb R$ then is $\Bbb R^2$ with the $x$-axis 
doubled: $(x,0)_+$ and $(x,0)_-$ cannot be separated for each $x\in \Bbb R$.
The charts $j_s(M)$ all are diffeomorphic to 
$M=\Bbb R^2\setminus\{0\}$ and contain $(x,0)_-$ for $x<0$ and 
$(x,0)_+$ for $x>0$. The charts $j_r(M)$ and $j_s(M)$ are glued 
together by the shift $x\mapsto x+s-r$. 
In this example $M_\Bbb R$ is not Hausdorff, but its Hausdorff 
quotient (given by the equivalence relation generated 
by identifying non-separable points) is again a smooth manifold and 
has the universal property described in \nmb!{1}. 

\subhead\nmb.{3}. Example \endsubhead
Let $(M,X)=(\Bbb R^2\setminus\{0\}\x [-1,1],\partial_x)$.
The trajectories of $\bar X$ on $\Bbb R\x M$ in the 
slices $y=\text{constant}$ for $|y|\le1$ and $|y|\ge1$ then look as 
in the second and third illustration above. The flow completion 
$M_\Bbb R$ then becomes $\Bbb R^2$ with the part $\Bbb R\x [-1,1]$ 
doubled and the topology such that the points $(x,-1)_-$ and 
$(x,-1)_+$ cannot be separated as well as the points $(x,1)_-$ and 
$(x,1)_+$. The flow is just $(x,y)\to (x+t,y)$:
$$
{\epsfxsize=8cm\epsfbox{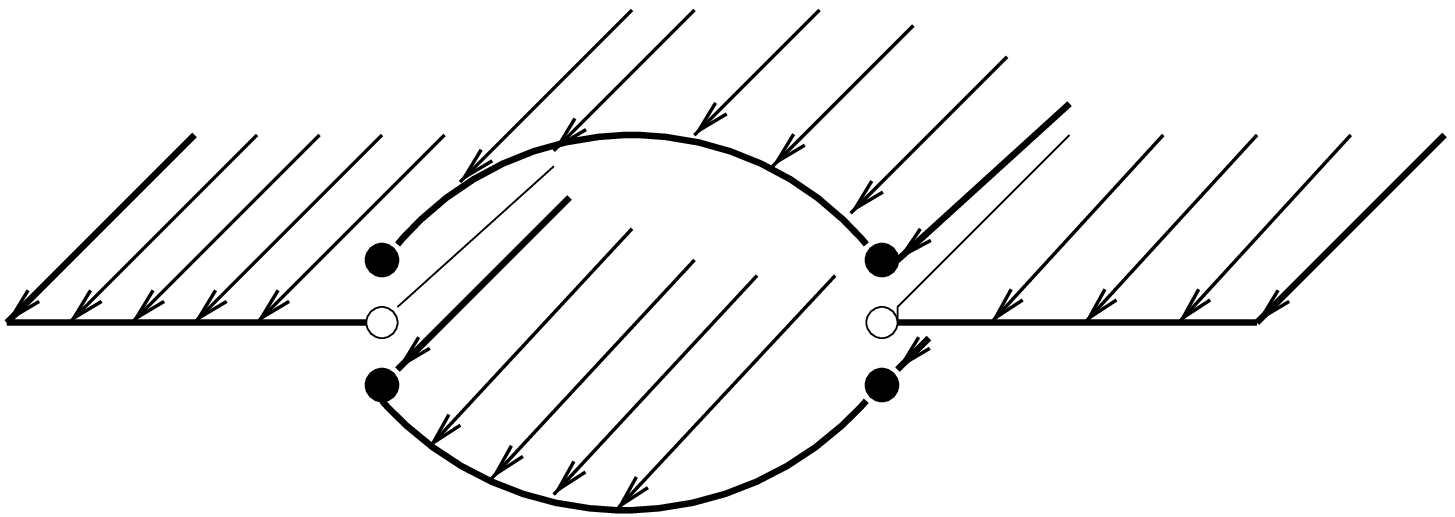}}
$$
In this example $M_\Bbb R$ is not Hausdorff, and its Hausdorff 
quotient is not a smooth manifold any more. There are two obvious 
quotient manifolds which are Hausdorff, the cylinder and the plane. 
Thus none of these two has the universal property of \nmb!{1}.

\subhead\nmb.{4}. Non-Hausdorff smooth manifolds \endsubhead
We met second countable smooth manifolds which need not be 
Hausdorff. Let us discuss a little their properties. They are $T_1$, 
since all points are closed; they are closed in a chart. The 
construction of the tangent bundle is by glueing the local tangent bundles. 
Smooth mappings and vector fields are defined as usual: Non separable 
pairs of points are mapped to non separable pairs. Vector fields 
admit flows as usual: These are given locally in the charts and are 
then glued together. If $x$ and $y$ are non separable points and if 
$X$ is a vector field on the manifold, then for each $t$ the points 
$\Fl^X_t(x)$ and $\Fl^X_t(y)$ are non separable. 
Theorem \nmb!{1} can be extended to the category of not necessarily 
Hausdorff smooth manifolds and vector fields, without any change in 
the proof. 

\subhead\nmb.{5}. Remark \endsubhead
The ideas in this paper generalize to the setting of 
$\frak g$-manifolds, where $\frak g$ is a finite dimension Lie group. 
Let $G$ be the simply connected Lie group with Lie algebra $\frak g$. 
Then one may construct the $G$-completion of a non-complete 
$\frak g$-manifold. There are difficulties with the property $T_1$, 
not only with Hausdorff.
This was our original road which was inspired by \cit!{1}. 
We treat the full theory in \cit!{2}.
We thought that the special case of a vector field is 
interesting in its own.

\enddocument